# Box-counting dimension and analytic properties of hidden variable fractal interpolation functions with function contractivity factors

Chol-Hui Yun[1)], Mi-Kyong Ri[2)]
[1), 2)] Faculty of Mathematics, Kim Il Sung University,
Pyongyang, Democratic People's Republic of Korea
*ch.yun@ryongnamsan.edu.kp

Abstract: We estimate the bounds of box-counting dimension of hidden variable fractal interpolation functions (HVFIFs) and hidden variable bivariate fractal interpolation functions (HVBFIFs) with four function contractivity factors and present analytic properties of HVFIFs which are constructed to ensure more flexibility and diversity in modeling natural phenomena. Firstly, we construct the HVFIFs and analyze their smoothness and stability. Secondly, we obtain the lower and upper bounds of box-counting dimension of the HVFIFs. Finally, in the similar way, we get the lower and upper bounds of box-counting dimension of HVBFIFs constructed in [21].

## 1. Introduction

Fractal interpolation function (FIF) is an interpolation function whose graph is a fractal set. In 1986, M. F. Barnsley [1] introduced a concept of FIF to model better natural phenomena which are irregular and complicated and the FIFs have been widely studied ever since in many papers. [2-15, 17-22]

In general, to get the FIF, we construct an iterated function system (IFS) on the basis of a given data set and then define a Read-Bajraktarevic operator on some space of continuous functions. A fixed point of the operator is an interpolation function of the data set and its graph is an attractor of the constructed IFS. The FIFs are not differentiable and fractal dimensions of their graphs are not integers. The fractal dimensions [10, 14, 15, 19 and 20], smoothness [5, 10, 11, 12 and 22] and stability [6, 7, 11, 13, 21 and 22] of the FIFs have been studied.

M.F. Barnsley *et al.* [2] introduced a concept of hidden variable fractal interpolation function (HVFIF) which is more complicated, diverse and irregular than the FIF for the same set of interpolation data. P. Bouboulis and L. Dalla [3] constructed hidden variable vector valued FIFs on grids in $\mathbf{R}^2$. Many researchers have studied the HVFIFs. [2-6, 12, 13, 15, 17, 18, 21 and 22]

The idea of the construction of the HVFIF is to extend the given data set on $\mathrm{R}^2$ into a data set on $\mathrm{R}^3$, make a vector valued fractal interpolation function for the extended data set and then project the vector valued function onto $\mathrm{R}^2$, which gives the HVFIF. It is usually non self-affine, because the HVFIF is the projection of a vector valued function. In [2], authors studied the HVFIF with three free parameters: free variable, constrained free variable, and free hidden variable (they are also called contractivity factors). Changing hidden variables, we can control shapes and fractal dimensions of the graphs of HVFIFs more flexibly.

A. K. B. Chand and G. P. Kapoor [5] studied a coalescence hidden variable fractal interpolation function using the IFS with free parameters and constrained free parameters. In many articles, authors have studied smoothness [5, 12 and 22], stability [6, 13, 21 and 22] and fractal dimension [5, 15] of HVFIFs.

The constructions of FIFs and HVFIFs with constant contractivity factors lack the flexibility which is necessary to model complicated and irregular natural phenomena.[1-7, 10, 12, 13, 15 and 18] A lot of fractal objects in nature have different contractivity factors at each point. Therefore, using function contractivity factors, we can model better mathematically fractal objects in nature. In many papers, the construction and analytic properties of FIFs [9, 11, 14, 19, and 20] and HVFIFs [16, 21, and 22] with function contractivity factors have been studied. The authors in [5, 15] estimated the lower and upper bounds of fractal dimension of HVFIF with three constant contractivity factors.

In order to ensure the flexibility and diversity of the methods of constructing FIFs, R. Uthayakumar and M. Rajkumar [17] presented a construction of HVFIF with function contractivity factors, where only one free variable is a function and the constrained free variable and free hidden variable are all constants. C.H. Yun and M.K. Li introduced a construction of hidden variable bivariate fractal interpolation functions (HVBFIFs) with four function contractivity factors [21] and studied their analytic properties. [22]

The main aim of the present paper is to find a lower and upper bounds of the fractal dimension of HVFIFs and HVBFIFs (constructed in [21]) with four function contractivity factors. To do this, first of all, we construct HVFIFs with four function contractivity factors and analyze their smoothness and stability, which are similar to ones in [21, 22]. Next, we estimate the box-counting dimension of the HVFIFs and HVBFIFs.

The remainder of this paper is organized as follows: In section 2, IFSs with function contractivity factors for the extended data set are constructed (Theorem 1) and using the Read-Bajraktarebic operator, vector valued fractal interpolation functions whose graphs are attractors of the IFSs are constructed. (Theorem 2) Then projecting them onto $\mathrm{R}^2$, HVFIFs for the given data set are constructed. In section 3, we consider smoothness (Theorem 3) and stability (Theorem 4-7) of the constructed HVFIFs. In section 4, the lower and upper bounds of fractal dimension of the HVFIFs are given (Theorem 8). In section 5, we introduce the result for the box-counting dimension of HVBFIFs constructed in [21].

## 2. Construction of hidden variable fractal interpolation function(HVFIF).

In this section, we construct hidden variable fractal interpolation functions with four function contractivity scaling factors.

Let a data set $P_0$ in $R^2$ be given by

$$P_0 = \{(x_i, y_i) \in \mathrm{R}^2; i = 0, 1, \cdots, n\}, \quad (-\infty < x_0 < x_1 < \cdots < x_n < +\infty).$$

To construct a hidden variable fractal interpolation function of this data set, we extend the data set $P_0$ into a data set $P$ in $\mathrm{R}^3$ as follows:

$$P=\{(x_i, y_i, z_i)=(x_i, \vec{y}_i)\in \mathrm{R}^3; i=0,1,\cdots,n\},\ (-\infty<x_0<x_1<\cdots<x_n<\infty),$$

where $\vec{y}_i=(y_i, z_i)$ and $z_i$ are parameters. Let us denote $N_n=\{1,2,\cdots,n\}$, $I_i=[x_{i-1}, x_i]$ and $I=[x_0, x_n]$.

Let mappings $L_i:[x_0, x_n]\to[x_{i-1}, x_i]$ for $i\in N_n$ be contraction homeomorphisms such that they map the end points of the interval $I$ into the end points of the subinterval $I_i$, i.e. $L_i(\{x_0, x_n\})=\{x_{i-1}, x_i\}$. We denote Lipschitz constant of Lipschitz mapping $f$ by $L_f(c_f)$.

We define mappings $\vec{F}_i: I\times \mathrm{R}^2\to \mathrm{R}^2$ for $i\in\{1,\cdots,n\}$ by

$$\vec{F}_i(x,\ \vec{y})=\begin{pmatrix} s_i(L_i(x))y+s_i'(L_i(x))z+q_i(x) \\ \tilde{s}_i(L_i(x))y+\tilde{s}_i'(L_i(x))z+\tilde{q}_i(x)\end{pmatrix}=\begin{pmatrix} s_i(L_i(x)) & s_i'(L_i(x)) \\ \tilde{s}_i(L_i(x)) & \tilde{s}_i'(L_i(x))\end{pmatrix}\begin{pmatrix} y \\ z\end{pmatrix}+\begin{pmatrix} q_i(x) \\ \tilde{q}_i(x)\end{pmatrix},$$

where functions $s_i, s_i', \tilde{s}_i, \tilde{s}_i': I_i\to R$ are arbitrary Lipschitz functions whose absolute values are less than 1 and $q_i$, $\tilde{q}_i: I\to R$ are Lipschitz functions defined as follows:

$$q_i(x)=-s_i(L_i(x))g_i(x)-s_i'(L_i(x))g_i'(x)+h_i(L_i(x)),$$
$$\tilde{q}_i(x)=-\tilde{s}_i(L_i(x))g_i(x)-\tilde{s}_i'(L_i(x))g_i'(x)+\tilde{h}_i(L_i(x)),$$

where $g_i$, $g_i': I\to R$ and $h_i$, $\tilde{h}_i: I_i\to R$ are Lipschitz functions such that

$$g_i(x_\alpha)=y_\alpha,\ \ g_i'(x_\alpha)=z_\alpha,\quad \alpha\in\{0, n\},$$
$$h_i(x_a)=y_a,\ \ \tilde{h}_i(x_a)=z_a,\quad a\in\{i-1, i\},$$

respectively.

For example, they can be constructed as the following Lagrange interpolation polynomials:

$$g_i(x)=\frac{x-x_0}{x_n-x_0}y_n+\frac{x-x_n}{x_0-x_n}y_0,\ \ g_i'(x)=\frac{x-x_0}{x_n-x_0}z_n+\frac{x-x_n}{x_0-x_n}z_0, \tag{1}$$
$$h_i(x)=\frac{x-x_{i-1}}{x_i-x_{i-1}}y_i+\frac{x-x_i}{x_{i-1}-x_i}y_{i-1},\ \ \tilde{h}_i(x)=\frac{x-x_{i-1}}{x_i-x_{i-1}}z_i+\frac{x-x_i}{x_{i-1}-x_i}z_{i-1}.$$

Then, we have

$$\vec{F}_i(x,\ \vec{y})=\begin{pmatrix} s_i(L_i(x))(y-g_i(x))+s_i'(L_i(x))(z-g_i'(x))+h_i(L_i(x)) \\ \tilde{s}_i(L_i(x))(y-g_i(x))+\tilde{s}_i'(L_i(x))(z-g_i'(x))+\tilde{h}_i(L_i(x))\end{pmatrix}$$

and for $\alpha\in\{0, n\}$ and $a\in\{i-1, i\}$ such that $L_i(x_\alpha)=x_a$, we get

$$\vec{F}_i(x_\alpha,\ \vec{y}_\alpha)=\vec{y}_a.$$

It is obvious that $\vec{F}_i(x,\ \vec{y})$ are Lipschitz mappings.

Let us denote as follows:

$$F_i^1(x, y, z)=s_i(L_i(x))y+s_i'(L_i(x))z+q_i(x),\ \ F_i^2(x, y, z)=\tilde{s}_i(L_i(x))y+\tilde{s}_i'(L_i(x))z+\tilde{q}_i(x).$$

Let us denote $\vec{F}_i(x,\ \vec{y})=\vec{S}_i\vec{y}+\vec{Q}_i(x)$, where

$$\vec{S}_i(x)=\begin{pmatrix} s_i(L_i(x)) & s_i'(L_i(x)) \\ \tilde{s}_i(L_i(x)) & \tilde{s}_i'(L_i(x)) \end{pmatrix},\ \vec{Q}_i(x)=\begin{pmatrix} q_i(x) \\ \tilde{q}_i(x) \end{pmatrix},\ \vec{y}=\begin{pmatrix} y \\ z \end{pmatrix}.$$

Suppose that $E \subset R^2$ is an enough large bounded domain containing $\vec{y}_i$, $i=1,\cdots,n$.

Now we define transformations $\vec{W}_i : I \times E \to I_i \times \mathrm{R}^2$, $i=1,\cdots,n$ by

$$\vec{W}_i(x,\vec{y})=(L_i(x),\ \vec{F}_i(x,\vec{y})),\ i=1,\cdots,n.$$

Then we can prove that $\vec{W}_i$ map the data points of $P$ on the end of interval $I$ into the data points of $P$ on the end of interval $I_i$. For a function *f*, let us denote $\|f\|_\infty=\sup_x |f(x)|$. We denote $S=\max\{\|s_i\|_\infty+\|\tilde{s}_i\|_\infty, \|s_i'\|_\infty+\|\tilde{s}_i'\|_\infty;\ i=1,\cdots,n\}$. The following theorem gives conditions for $\vec{W}_i$ to be contraction transformations.

**Theorem 1.** *If* $S<1$, *then there exists a distance* $\rho_\theta$ *equivalent to the Euclidean metric on* $\mathrm{R}^3$ *such that* $\vec{W}_i$ *are contraction mappings with respect to the distance* $\rho_\theta$.

**Proof.** We take $\theta$ as a positive number satisfying

$$\theta<\frac{1-c_L}{L_S\kappa+L_Q},$$

where $c_L=\max\{c_{L_i};i=1,\cdots,n\}$, $L_S=\max\{L_{s_i}c_{L_i}+L_{\tilde{s}_i}c_{L_i}, L_{s_i'}c_{L_i}+L_{\tilde{s}_i'}c_{L_i};i=1,\cdots,n\}$, $\kappa=\sup_{\vec{y}\in D}\|\vec{y}\|_1$, $L_Q=\max\{L_{q_i}+L_{\tilde{q}_i};i=1,\cdots,n\}$.

$\|\cdot\|_1$ is a norm on $R^2$. Let us define a distance $\rho_\theta$ on $\mathrm{R}^3$ by

$$\rho_\theta((x,\vec{y}),(x',\vec{y}'))=\|x-x'\|_1+\theta\|\vec{y}-\vec{y}'\|_1,\quad (x,\vec{y}),(x',\vec{y}')\in\mathrm{R}^3.$$

It is clear that the $\rho_\theta$ is equivalent to Euclidean metric on $\mathrm{R}^3$. For $(x,\vec{y}),(x',\vec{y}')\in I\times D$, we have

$$\begin{aligned}
&\rho_\theta(\vec{W}_i(x,\vec{y}),\vec{W}_i(x',\vec{y}'))=\|\vec{L}_i(x)-\vec{L}_i(x')\|_1+\theta\|\vec{F}_i(x,\vec{y})-\vec{F}_i(x',\vec{y}')\|_1\\
&=\|L_i(x)-L_i(x')\|_1+\theta\|\vec{F}_i(x,\vec{y})-\vec{F}_i(x',\vec{y}')\|_1\\
&=\|L_i(x)-L_i(x')\|_1+\theta\|\vec{S}_i(x)\vec{y}+\vec{Q}_i(x)-\vec{S}_i(x')\vec{y}'-\vec{Q}_i(x')\|_1\\
&=\|L_i(x)-L_i(x')\|_1+\theta\|(\vec{S}_i(x)\vec{y}-\vec{S}_i(x)\vec{y}')+(\vec{S}_i(x)\vec{y}'-\vec{S}_i(x')\vec{y}')+(\vec{Q}_i(x)-\vec{Q}_i(x'))\|_1\\
&\le\|L_i(x)-L_i(x')\|_1+\theta(\|\vec{S}_i(x)\|_1\|\vec{y}-\vec{y}'\|_1+\|\vec{S}_i(x)-\vec{S}_i(x')\|_1\|\vec{y}'\|_1+\|\vec{Q}_i(x)-\vec{Q}_i(x')\|_1)\\
&\le c_L\|x-x'\|_1+S\theta\|\vec{y}-\vec{y}'\|_1+L_S\kappa\theta\|x-x'\|_1+L_Q\theta\|x-x'\|_1\\
&=(c_L+\theta(L_S\kappa+L_Q))\|x-x'\|_1+S\theta\|\vec{y}-\vec{y}'\|_1\\
&\le\max\{c_L+\theta(L_S\kappa+L_Q),S\}(\|x-x'\|_1+\theta\|\vec{y}-\vec{y}'\|_1)\\
&=c\rho_\theta((x,\vec{y}),(x',\vec{y}')).
\end{aligned}$$

From the hypothesis of the theorem and the condition on $\theta$, it follows that $c=\max\{c_L+\theta$ $(L_S\alpha+L_Q),\ S\}<1$. Therefore, $\vec{W}_i$ are contraction transformations. □

*Remark 1.* Even in the case when $\|\cdot\|_1$ in the definition of $\rho_\theta$ is changed into $\|\cdot\|_\infty$, we get the similar result to Theorem 1.

Therefore, $\{R^3; W_i, i=1,\cdots,n\}$ is a hyperbolic iterated function system (IFS) corresponding to the extended data set $P$. Let us denote by $A$ an attractor of the IFS.

For the IFS, we have the following theorem.

**Theorem 2.** *There is a continuous interpolation function* $\vec{f}$ *of the extended data set* $P$ *such that the graph of* $\vec{f}$ *is the attractor of the above constructed IFS.*

**Proof.** Let us define a set $\overline{C}(I)$ as follows:

$$\overline{C}(I)=\{\vec{h}:I\to R^2;\vec{h}\ \text{ interpolates the data set } P \text{ and is continuous.}\}.$$

Then, it can be easily proved that the set $\overline{C}(I)$ is a complete metric space with respected to a norm $\|\cdot\|_\infty$.

For $\vec{h}(\in\overline{C}(I))$, we define a mapping $T\vec{h}$ on $I$ as follows:

$$(T\vec{h})(x)=\vec{F}_i(L_i^{-1}(x),\ \vec{h}(L_i^{-1}(x))),\quad x\in I_i.$$

Then, it follows that $T\vec{h}\in\overline{C}(I)$.

In fact, for any $i\in\{0,1,\cdots,n\}$, there is $\alpha\in\{0,n\}$ such that

$$(T\vec{h})(x_i)=\vec{F}_i(L_i^{-1}(x_i),\ \vec{h}(L_i^{-1}(x_i)))=\vec{F}_i(x_\alpha,\vec{y}_\alpha)=\vec{y}_i,$$

where $L_i(x_\alpha)=x_i$ and $L_0(x)=L_1(x)$.

Therefore, an operator $T:\overline{C}(E)\to\overline{C}(E)$ is well defined on $\overline{C}(E)$. Moreover, it is a contraction operator, because

$$\|(T\vec{h})(x)-(T\vec{h}')(x)\|_1=\|\vec{F}_i(L_i^{-1}(x),\vec{h}(L_i^{-1}(x)))-\vec{F}_i(L_i^{-1}(x),\vec{h}'(L_i^{-1}(x)))\|_1$$
$$=\|\vec{S}_i(x)\vec{h}(L_i^{-1}(x))+\vec{Q}_i(L_i^{-1}(x))-\vec{S}_i(x)\vec{h}'(L_i^{-1}(x))-\vec{Q}_i(L_i^{-1}(x))\|_1$$
$$=\|\vec{S}_i(x)\|_1\|\vec{h}(L_i^{-1}(x))-\vec{h}'(L_i^{-1}(x))\|_1$$
$$\le S\|\vec{h}-\vec{h}'\|_1.$$

Hence, from the fixed point theorem in the complete metric space, it follows that $T$ has a unique fixed point $\vec{f}(\in\overline{C}(I))$. Then, we have

$$\vec{f}(x)=\vec{F}_i(L_i^{-1}(x),\ \vec{f}(L_i^{-1}(x))).$$

Therefore, for the graph $Gr(\vec{f})$ of $\vec{f}$, we have

$$Gr(\vec{f})=\bigcup_{j=1}^{n}W_j(Gr(\vec{f})).$$

This means that $Gr(\vec{f})$ is the attractor of the constructed IFS. Therefore, from the uniqueness of an attractor, we have $A=Gr(\vec{f})$.□

Let us denote the vector valued function $\vec{f}:I\to R^2$ by $\vec{f}=(f_1(x),\ f_2(x))$, where $f_1:I\to R$ interpolates the given data set $P_0$ and is called a *hidden variable fractal interpolation function* (HVFIF) of the data set $P_0$. Then, a set $\{(x,\ f_1(x)): x \in I\}$ is a projection of $A$ onto $R^2$. Since a projection of the attractor is not always self-affine set, the hidden variable fractal interpolation function is not self-affine fractal interpolation function in general. The function $f_2(x)$ interpolates a set $\{(x_i, z_i,); (x_i, y_i, z_i) \in P,\ \ i = 0, 1, \cdots, n\}$. As we can know from the above proof, we get

$$\vec{f}(x, y) = \vec{F}_i(L_i^{-1}(x),\ \vec{f}(L_i^{-1}(x))),\quad x \in I_i,$$

i.e.

$$\vec{f}(x, y) = \vec{F}_i(L_i^{-1}(x),\ f_1(L_i^{-1}(x)),\ f_2(L_i^{-1}(x))),\quad x \in I_i.$$

Therefore, for all $x \in I$, we have

$$f_1(x) = s_i(x) f_1(L_i^{-1}(x)) + s_i'(x) f_2(L_i^{-1}(x)) + q_i(L_i^{-1}(x)), \tag{2}$$

$$f_2(x) = \tilde{s}_i(x) f_1(L_i^{-1}(x)) + \tilde{s}_i'(x) f_2(L_i^{-1}(x)) + \tilde{q}_i(L_i^{-1}(x)). \tag{3}$$

*Remark 2.* The contractivity factors $s_i$, $s_i'$, $\tilde{s}_i'$ of the above IFS are constants and $\tilde{s}_i = 0$ in [3- 6, 12, 13]. In [17], $s_i$ are functions and $s_i'$, $\tilde{s}_i'$ are constants and $\tilde{s}_i = 0$.

Example. Fig.1 shows the graphs of HVFIFs constructed from a data set $P$ ={(1, 20), (0.25, 30), (0.5, 10), (0.75, 50), (1, 40)} with different contractivity factor functions. The contractivity factors $\{s_1,\ s_2,\ s_3, s_4\}$, $\{\tilde{s}_1, \tilde{s}_2, \tilde{s}_3, \tilde{s}_4\}$, $\{s_1',\ s_2',\ s_3', s_{41}'\}$, $\{\tilde{s}_1', \tilde{s}_2', \tilde{s}_3', \tilde{s}_4'\}$ are as follows: (a) {0.3, 0.85, 0.8, 0.5}, {0, 0, 0, 0}, {0.8, 0.6, 0.4, 0.5 }, {0.19, 0.37, 0.48, 0.43}, (b) {0.3, 0.85, 0.8, 0.5}, {0.64, 0.14, 0.19, 0.49}, {0.8, 0.6, 0.4, 0.5 }, {0.19, 0.37, 0.48, 0.43}, (c) {sin($x$), cos(30 $x$), sin($x$), cos(5 $x$)}, {0, 0, 0, 0}, {2.9x, 1.9$x$, $x$, $x$ }, {0.9-2.9 $x$, 0.95-1.9 $x$, 0.9- $x$, 0.99- $x$ }, (d) {sin($x$), cos(30 $x$), sin($x$), cos(5 $x$)}, {0.9-| sin($x$)|, 0.89-| cos(30$x$)|, 0.94-| sin($x$)|, 0.9-| cos(5 $x$)|}, {2.9 $x$, 1.9 $x$, $x$, $x$ }, {0.9-2.9 $x$, 0.95-1.9 $x$, 0.9- $x$, 0.99-$x$}.

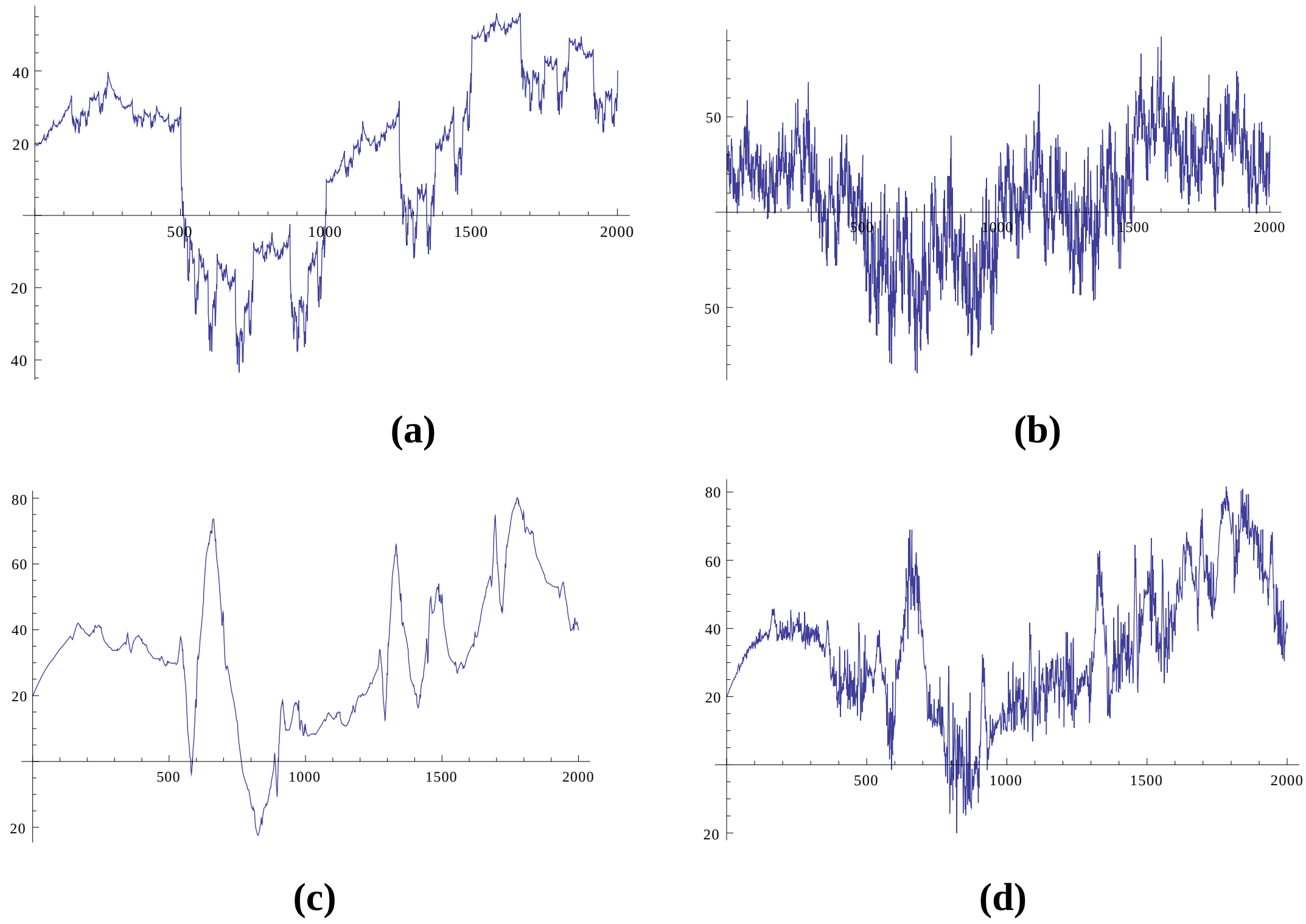

Fig.1. Hidden variable fractal interpolation functions.

## 3. Smoothness and stability of the HVFIFs.

In this section, we analyze the smoothness and stability of the constructed HVFIF.

Let $I=[0,1]$ and $L_i : I \to I_i$ be similarities with Lipschitz constants $L_{L_i}$. Let us denote as follows: $\bar{g}=\max_x |g(x)|$ , $L_s=\max_i\{L_{s_i}\}$, $L_{s'}=\max_i\{L_{s_i'}\}$, $L_{\tilde{s}}=\max_i\{L_{s_i}\}$, $L_{\tilde{s}'}=\max_i\{L_{s_i'}\}$,

$L_q=\max_i\{L_{q_i}\}, L_{\tilde{q}}=\max_i\{L_{\tilde{q}_i}\}$, $|I_{\min}|=\min_i\{|I_i|\}$, $|I_{\max}|=\max_i\{|I_i|\}$,

$\bar{\omega}_k=\max_{x\in I}\{|s_k(x)|, |s_k'(x)|\}$, $\bar{\tilde{\omega}}_k=\max_{x\in I}\{|\tilde{s}_k(x)|, |\tilde{s}_k'(x)|\}$, $\omega=\max_k\{\bar{\omega}_k\}$, $\tilde{\omega}=\max_k\{\bar{\tilde{\omega}}_k\}$.

**Lemma 1**. ([21]) *If* $0<\alpha$, *then for* $x$ obeying $0<x<1$,

$$0<-x^{\alpha}\ln x\le \frac{1}{\alpha e}.$$

**Theorem 3.** *Let* $f_1(x)$ *be the HVFIF constructed in the Theorem 2 and* $\max\{\omega,\ \tilde{\omega}\}<\frac{|I_{\min}|}{2|I_{\max}|}$. *Then there exist a positive constant* $L_1$ *and* $\tau_1$ *with* $0<\tau_1\le 1$ *such that*

$$|f_1(x)-f_1(\bar{x})|\le L_1|x-\bar{x}|^{\tau_1}. \tag{4}$$

**Proof.** Let us denote $I_{r_1r_2\dots r_m}=L_{r_m}\circ L_{r_{m-1}}\circ\cdots\circ L_{r_1}(I)$, $I_0=I$. Then, for any $x$, $\bar{x}\in[0,1]$ with $x\le\bar{x}$, there is $m\in N$ such that $I_{r_1r_2\dots r_m}\subset[x,\bar{x}]\subset I_{r_k\dots r_m}$, $r_1,r_2,\dots,r_m\in N_n$, $k=1,2$ and we get

$$\begin{aligned}
|f_1(x)-f_1(\bar{x})|&=|f_1(L_{r_m}(L_{r_m}^{-1}(x)))-f_1(L_{r_m}(L_{r_m}^{-1}(\bar{x})))|\\
&=|(s_{r_m}(x)f_1(L_{r_m}^{-1}(x))+s'_{r_m}(x)f_2(L_{r_m}^{-1}(x))+q_{r_m}(L_{r_m}^{-1}(x)))-(s_{r_m}(\bar{x})f_1(L_{r_m}^{-1}(\bar{x}))\\
&\qquad+s'_{r_m}(\bar{x})f_2(L_{r_m}^{-1}(\bar{x}))+q_{r_m}(L_{r_m}^{-1}(\bar{x})))|\\
&\le|s_{r_m}(x)f_1(L_{r_m}^{-1}(x))-s_{r_m}(\bar{x})f_1(L_{r_m}^{-1}(\bar{x}))|+|s'_{r_m}(x)f_2(L_{r_m}^{-1}(x))-s'_{r_m}(\bar{x})f_2(L_{r_m}^{-1}(\bar{x}))|\\
&\qquad+|q_{r_m}(L_{r_m}^{-1}(x))-q_{r_m}(L_{r_m}^{-1}(\bar{x}))|,
\end{aligned} \tag{5}$$

and

$$\begin{aligned}
&|s_{r_m}(x)f_1(L_{r_m}^{-1}(x))-s_{r_m}(\bar{x})f_1(L_{r_m}^{-1}(\bar{x}))|\\
&\le|s_{r_m}(x)f_1(L_{r_m}^{-1}(x))-s_{r_m}(x)f_1(L_{r_m}^{-1}(\bar{x}))|+|s_{r_m}(x)f_1(L_{r_m}^{-1}(\bar{x}))-s_{r_m}(\bar{x})f_1(L_{r_m}^{-1}(\bar{x}))|\\
&\le|s_{r_m}(x)|\cdot|f_1(L_{r_m}^{-1}(x))-f_1(L_{r_m}^{-1}(\bar{x}))|+|s_{r_m}(x)-s_{r_m}(\bar{x})|\cdot|f_1(L_{r_m}^{-1}(\bar{x}))|\\
&\le|s_{r_m}(x)|\cdot|f_1(L_{r_m}^{-1}(x))-f_1(L_{r_m}^{-1}(\bar{x}))|+L_{s_{r_m}}|x-\bar{x}|\cdot\|f_1\|_\infty
\end{aligned} \tag{6}$$

$$\begin{aligned}
&|s'_{r_m}(x)f_2(L_{r_m}^{-1}(x))-s'_{r_m}(\bar{x})f_2(L_{r_m}^{-1}(\bar{x}))|\\
&\le|s'_{r_m}(x)f_2(L_{r_m}^{-1}(x))-s'_{r_m}(x)f_2(L_{r_m}^{-1}(\bar{x}))|+|s'_{r_m}(x)f_2(L_{r_m}^{-1}(\bar{x}))-s'_{r_m}(\bar{x})f_2(L_{r_m}^{-1}(\bar{x}))|\\
&\le|s'_{r_m}(x)|\cdot|f_2(L_{r_m}^{-1}(x))-f_2(L_{r_m}^{-1}(\bar{x}))|+|s'_{r_m}(x)-s'_{r_m}(\bar{x})|\cdot|f_2(L_{r_m}^{-1}(\bar{x}))|\\
&\le|s'_{r_m}(x)|\cdot|f_2(L_{r_m}^{-1}(x))-f_2(L_{r_m}^{-1}(\bar{x}))|+L_{s'_{r_m}}|x-\bar{x}|\cdot\|f_2\|_\infty
\end{aligned} \tag{7}$$

$$|q_{r_m}(L_{r_m}^{-1}(x))-q_{r_m}(L_{r_m}^{-1}(\bar{x}))|\le L_{q_{r_m}}|L_{r_m}^{-1}(x)-L_{r_m}^{-1}(\bar{x})|=\frac{L_{q_{r_m}}}{|I_{r_m}|}|x-\bar{x}|. \tag{8}$$

Therefore, by (6), (7) and (8), it follows that

$$|f_1(x)-f_1(\bar{x})|\le\left(L_{s_{r_m}}\|f_1\|_\infty+L_{s'_{r_m}}\|f_2\|_\infty+\frac{L_{q_{r_m}}}{|I_{r_m}|}\right)|x-\bar{x}|+$$

$$+|s_{r_m}(x)|\cdot|f_1(L_{r_m}^{-1}(x))-f_1(L_{r_m}^{-1}(\bar{x}))|+|s'_{r_m}(x)|\cdot|f_2(L_{r_m}^{-1}(x))-f_2(L_{r_m}^{-1}(\bar{x}))|=$$

$$=M_{r_m}|x-\bar{x}|+\bar{\omega}_{r_m}(|f_1(L_{r_m}^{-1}(x))-f_1(L_{r_m}^{-1}(\bar{x}))|+|f_2(L_{r_m}^{-1}(x))-f_2(L_{r_m}^{-1}(\bar{x}))|),$$

where

$$M_k=L_{s_k}\|f_1\|_\infty+L_{s'_k}\|f_2\|_\infty+\frac{L_{q_k}}{|I_k|}. \tag{9}$$

Since $I_{r_1r_2\dots r_m}\subset[x,\bar{x}]\subset I_{r_k\dots r_m}$, we have

$$|I_{\max}|^{m-2}\ge|I_{r_k}|\cdots|I_{r_m}|\ge|x-\bar{x}|\ge|I_{r_1}|\cdot|I_{r_2}|\cdots|I_{r_m}|\ge|I_{\min}|^m. \tag{10}$$

Let us denote as follows:

$$L_{r_1r_2\cdots r_m}^{-1}:=L_{r_1}^{-1}\circ L_{r_2}^{-1}\circ\cdots\circ L_{r_m}^{-1},\quad \delta:=\max_k\left\{\frac{2\max\{\bar{\omega}_k,\tilde{\bar{\omega}}_k\}}{|I_k|}\right\},$$

$$\tilde{M}_k:=L_{\tilde{s}_k}\|f_1\|_\infty+L_{\tilde{s}'_k}\|f_2\|_\infty+\frac{L_{\tilde{q}_k}}{|I_k|},\quad M:=\max_k\{M_k,\tilde{M}_k\}. \tag{11}$$

Then, by induction, we have

$$|f_1(L_{r_m}^{-1}(x))-f_1(L_{r_m}^{-1}(\bar{x}))|=|f_1(L_{r_{m-1}}^{-1}(L_{r_{m-1}r_m}^{-1}(x)))-f_1(L_{r_m}^{-1}(L_{r_{m-1}r_m}^{-1}(\bar{x})))|$$

$$\le M_{r_{m-1}}\left|L_{r_m}^{-1}(x)-L_{r_m}^{-1}(\bar{x})\right|+\bar{\omega}_{r_{m-1}}(|f_1(L_{r_{m-1}r_m}^{-1}(x))-f_1(L_{r_{m-1}r_m}^{-1}(\bar{x}))|+|f_2(L_{r_{m-1}r_m}^{-1}(x))-f_2(L_{r_{m-1}r_m}^{-1}(\bar{x}))|)$$

$$\le\frac{M_{r_{m-1}}}{|I_{r_m}|}|x-\bar{x}|+\bar{\omega}_{r_{m-1}}(|f_1(L_{r_{m-1}r_m}^{-1}(x))-f_1(L_{r_{m-1}r_m}^{-1}(\bar{x}))|+|f_2(L_{r_{m-1}r_m}^{-1}(x))-f_2(L_{r_{m-1}r_m}^{-1}(\bar{x}))|).$$

Similarly, we get

$$|f_2(L_{r_m}^{-1}(x))-f_2(L_{r_m}^{-1}(\bar{x}))|\le$$

$$\le\frac{\tilde{M}_{r_{m-1}}}{|I_{r_m}|}|x-\bar{x}|+\tilde{\bar{\omega}}_{r_{m-1}}(|f_1(L_{r_{m-1}r_m}^{-1}(x))-f_1(L_{r_{m-1}r_m}^{-1}(\bar{x}))|+|f_2(L_{r_{m-1}r_m}^{-1}(x))-f_2(L_{r_{m-1}r_m}^{-1}(\bar{x}))|).$$

Then,

$$|f_1(x)-f_1(\bar{x})|\le$$

$$\le M_{r_m}|x-\bar{x}|+\bar{\omega}_{r_m}\left(\frac{M_{r_{m-1}}}{|I_{r_m}|}|x-\bar{x}|+\bar{\omega}_{r_{m-1}}(|f_1(L_{r_{m-1}r_m}^{-1}(x))-f_1(L_{r_{m-1}r_m}^{-1}(\bar{x}))|+|f_2(L_{r_{m-1}r_m}^{-1}(x))-f_2(L_{r_{m-1}r_m}^{-1}(\bar{x}))|)+\right.$$

$$\left.+\frac{\tilde{M}_{r_{m-1}}}{|I_{r_m}|}|x-\bar{x}|+\tilde{\bar{\omega}}_{r_{m-1}}(|f_1(L_{r_{m-1}r_m}^{-1}(x))-f_1(L_{r_{m-1}r_m}^{-1}(\bar{x}))|+|f_2(L_{r_{m-1}r_m}^{-1}(x))-f_2(L_{r_{m-1}r_m}^{-1}(\bar{x}))|)\right)\le$$

$$\le M|x-\bar{x}|+\frac{2\bar{\omega}_{r_m}}{|I_{r_m}|}M|x-\bar{x}|+\bar{\omega}_{r_m}(\bar{\omega}_{r_{m-1}}+\tilde{\bar{\omega}}_{r_{m-1}})(|f_1(L_{r_{m-1}r_m}^{-1}(x))-f_1(L_{r_{m-1}r_m}^{-1}(\bar{x}))|+$$

$$+|f_2(L_{r_{m-1}r_m}^{-1}(x))-f_2(L_{r_{m-1}r_m}^{-1}(\bar{x}))|)$$

$$\le\left(M+\frac{2\bar{\omega}_{r_m}}{|I_{r_m}|}M+\frac{2\bar{\omega}_{r_m}}{|I_{r_m}|}\frac{\bar{\omega}_{r_{m-1}}+\tilde{\bar{\omega}}_{r_{m-1}}}{|I_{r_{m-1}}|}M\right)|x-\bar{x}|+\bar{\omega}_{r_m}(\bar{\omega}_{r_{m-1}}+\tilde{\bar{\omega}}_{r_{m-1}})(\bar{\omega}_{r_{m-2}}+\tilde{\bar{\omega}}_{r_{m-2}})$$

$$(|f_1(L_{r_{m-2}r_{m-1}r_m}^{-1}(x))-f_1(L_{r_{m-2}r_{m-1}r_m}^{-1}(\bar{x}))|+|f_2(L_{r_{m-2}r_{m-1}r_m}^{-1}(x))-f_2(L_{r_{m-2}r_{m-1}r_m}^{-1}(\bar{x}))|$$

$$\le\cdots\le$$

$$\le \left( M + \frac{2\overline{\omega}_{r_m}}{|I_{r_m}|} M + \frac{2\overline{\omega}_{r_m}}{|I_{r_m}|} \frac{\overline{\omega}_{r_{m-1}} + \overline{\overline{\omega}}_{r_{m-1}}}{|I_{r_{m-1}}|} M + \cdots + \left( \frac{2\overline{\omega}_{r_m}}{\overline{\omega}_{r_m} + \overline{\overline{\omega}}_{r_m}} \prod_{k=4}^{m} \frac{\overline{\omega}_{r_k} + \overline{\overline{\omega}}_{r_k}}{|I_{r_k}|} \right) M \right) |x - \overline{x}| +$$

$$+ \left( \prod_{k=3}^{m-1} (\overline{\omega}_{r_k} + \overline{\overline{\omega}}_{r_k}) \right) \overline{\omega}_{r_m} (| f_1(L_{r_3\cdots r_m}^{-1}(x)) - f_1(L_{r_3\cdots r_m}^{-1}(\overline{x}))| + | f_2(L_{r_3\cdots r_m}^{-1}(x)) - f_2(L_{r_3\cdots r_m}^{-1}(\overline{x}))|)$$

$$\le M(1 + \delta + \delta^2 + \cdots + \delta^{m-3}) |x - \overline{x}| + \left( \prod_{k=3}^{m-1} (\overline{\omega}_{r_k} + \overline{\overline{\omega}}_{r_k}) \right) \overline{\omega}_{r_m} (2 \| f_1 \|_\infty + 2 \| f_2 \|_\infty)$$

$$\le M(1 + \delta + \delta^2 + \cdots + \delta^{m-3}) |x - \overline{x}| + 2(\| f_1 \|_\infty + \| f_2 \|_\infty) \left( \prod_{k=3}^{m} (\overline{\omega}_{r_k} + \overline{\overline{\omega}}_{r_k}) \right)$$

$$\le M(1 + \delta + \delta^2 + \cdots + \delta^{m-3}) |x - \overline{x}| + 2(\| f_1 \|_\infty + \| f_2 \|_\infty) \left( \prod_{k=3}^{m} |I_{r_k}| \right) \delta^{m-2}$$

$$\le M(1 + \delta + \delta^2 + \cdots + \delta^{m-3}) |x - \overline{x}| + \frac{2(\| f_1 \|_\infty + \| f_2 \|_\infty)}{|I_{r_1}| \cdot |I_{r_2}|} \delta^{m-2} |x - \overline{x}|$$

$$\le D \left( \sum_{k=0}^{m-2} \delta \right) |x - \overline{x}|, \tag{12}$$

where $D = \max \left\{ M, \ \frac{2(\| f_1 \|_\infty + \| f_2 \|_\infty)}{l_L^2} \right\}$.

(1) If $\delta < 1$, then by (12), we have

$$D \left( \sum_{k=0}^{m-2} \delta \right) |x - \overline{x}| = D \frac{1 - \delta^{m-1}}{1 - \delta} |x - \overline{x}| < D \frac{1}{1 - \delta} |x - \overline{x}|.$$

Hence, $| f_1(x) - f_1(\overline{x})| \le D \frac{1}{1 - \delta} |x - \overline{x}|$. Therefore, let us denote $L_1 = D \frac{1}{1 - \delta}$, $\tau_1 = 1$, then

$$| f_1(x) - f_1(\overline{x})| \le L_1 |x - \overline{x}|^{\tau_1}.$$

(2) If $\delta = 1$, then $D \left( \sum_{k=0}^{m-2} \delta \right) = D(m-1)$ by (12) and $|x - \overline{x}| \le |I_{\max}|^{m-2} < 1$ by (10). Therefore,

$$m - 2 \le \frac{\ln |x - \overline{x}|}{\ln |I_{\max}|}.$$

This gives

$$| f_1(x) - f_1(\overline{x})| \le D(m-1) |x - \overline{x}| \le D \left( 1 + \frac{\ln |x - \overline{x}|}{\ln |I_{\max}|} \right) |x - \overline{x}|$$

$$= D \left[ |x - \overline{x}| + \left( - \frac{|x - \overline{x}|^\alpha \ln |x - \overline{x}|}{|\ln |I_{\max}||} \right) |x - \overline{x}|^{1-\alpha} \right],$$

where $0 < \alpha < 1$. Since $0 < |x - \overline{x}| < 1$, by Lemma1, we have

$$| f_1(x) - f_1(\overline{x})| \le D \left( |x - \overline{x}| + \frac{D}{\alpha e |\ln |I_{\max}||} |x - \overline{x}|^{1-\alpha} \right).$$

Let us denote $L_1 = D \left( 1 + \frac{1}{\alpha e |\ln |I_{\max}||} \right)$, $\tau_1 = 1 - \alpha < 1$. Then

$$| f_1(x) - f_1(\overline{x})| \le L_1 |x - \overline{x}|^{\tau_1}.$$

(3) If $\delta > 1$, then

$$D\left(\sum_{k=0}^{m-2}\delta\right) = D\frac{\delta^{m-2}\left(1-\frac{1}{\delta^{m-1}}\right)}{1-\frac{1}{\delta}} < D\frac{\delta^{m-2}}{1-\frac{1}{\delta}} = D\frac{\delta^{m-1}}{\delta-1} \quad . \tag{13}$$

Let us choose $\tau$ such that $0 < \tau \le \frac{\ln\delta}{\ln L_L} + 1$. Then

$$\delta^m |x-\bar{x}| \le |x-\bar{x}|^\tau . \tag{14}$$

In fact, we get

$$\frac{\ln(\delta^m |x-\bar{x}|)}{\ln|x-\bar{x}|} = \frac{m\ln\delta}{\ln|x-\bar{x}|} + 1$$

and by (10), it follows that $\ln|x-\bar{x}| \le (m-2)\ln L_L$. Hence,

$$\frac{m\ln\delta}{\ln|x-\bar{x}|} + 1 \ge \frac{m\ln\delta}{(m-2)\ln|I_{\max}|} + 1 \ge \frac{\ln\delta}{\ln|I_{\max}|} + 1,$$

$$\frac{\ln(\delta^m |x-\bar{x}|)}{\ln|x-\bar{x}|} \ge \frac{\ln\delta}{\ln|I_{\max}|} + 1 .$$

Therefore, by (12), (13) and (14), choosing $L_1 = D/(\delta-1)$ and $\tau_1 = \tau$ gives

$$|f_1(x) - f_1(\bar{x})| \le L_1 |x-\bar{x}|^{\tau_1} . \quad \square$$

*Remark 3*. For the interpolation function $f_2(x)$, we can prove similarly that (4) holds true, i.e. there exist positive number $L_2$ and $0 < \tau_2 < 1$ such that for any $x,\ \bar{x} \in I$,

$$|f_2(x) - f_2(\bar{x})| \le L_2 |x-\bar{x}|^{\tau_2}$$

Next let us consider the stability of the HVFIF. Firstly, we mention the stability of the HVFIF according to perturbation of each coordinate of points in the data set *P*, and then all coordinates.

Let $I = [0,1]$. For $x_0 = x_0^* < x_1^* < \cdots < x_n^* = x_n$, let us denote $I_i^* = [x_{i-1}^*, x_i^*]$ and define mapping $R: I \to I$ as follows:

$$R(x) = x_{i-1}^* + \frac{x_i^* - x_{i-1}^*}{x_i - x_{i-1}}(x - x_{i-1}), \text{ for } x \in I_i .$$

Then $R(I_i) = I_i^*$.

Now, let us denote as follows:

$$L_i^* = R \circ L_i,\ s_i^* = s_i \circ R^{-1},\ s_i'^* = s_i' \circ R^{-1},\ \tilde{s}_i^* = \tilde{s}_i \circ R^{-1},$$
$$\tilde{s}_i'^* = \tilde{s}_i' \circ R^{-1},\ q_i^* = q_i \circ R^{-1},\ \tilde{q}_i^* = \tilde{q}_i \circ R^{-1}. \quad .$$

Let us denote by $f_1^*$ the HVFIF constructed using the construction in section 2 with (1) and an extended data

$$P_{x^*} = \{(x_i^*,\ y_i,\ z_i) \in \mathrm{R}^3; i = 0,1,\cdots,n\}(x_0 = x_0^* < x_1^* < \cdots < x_n^* = x_n) .$$

**Theorem 4**. *Let* $f_1,\ f_1^*$ *be HFIF for the extended data sets* $P$ *and* $P_{x^*}$, *respectively, and*

$\max\{\omega,\ \tilde{\omega}\} < \dfrac{|I_{\min}|}{2|I_{\max}|}$. *Then*

$$\| f_1 - f_1^* \|_\infty \le \frac{(1-\tilde{\omega})L_1 + \omega L_2 + (1-\tilde{\omega})L_q + \omega L_{\tilde{q}}}{1-\omega-\tilde{\omega}} \max_i |x_i - x_i^*|^\tau,\ \ \tau = \max\{\tau_1,\ \tau_2\},$$

*where* $L_1$ *and* $\tau_1$ *are constants in Theorem 3 and* $L_2$ *and* $\tau_2$ *are constants in Remark 3.*

We abbreviate the proof of Theorem 5. This can be proved in the similar way of Theorem 3 in [21].

Let us denote by $f_1^*$ the HVFIF constructed using the construction in section 2 for an extended data set $P_{y^*} = \{(x_i,\ y_i^*,\ z_i) \in \mathrm{R}^3;\ i = 0, 1, \cdots, n\}$.

**Theorem 5.** *Let* $f_1,\ f_1^*$ *be HFIF for the extended data sets* $P$ *and* $P_{y^*}$, *respectively, and* $\max\{\omega,\ \tilde{\omega}\} < \dfrac{|I_{\min}|}{2|I_{\max}|}$. *Then*

$$\| f_1 - f_1^* \|_\infty \le \frac{1+2\omega-\tilde{\omega}}{1-\omega-\tilde{\omega}} \max_i |y_i - y_i^*|.$$

Let us denote by $f_1^*$ the HVFIF constructed using the construction in section 2 for an extended data set $P_{z^*} = \{(x_i,\ y_i,\ z_i^*) \in \mathrm{R}^3;\ i = 0, 1, \cdots, n\}$.Then, the following theorem holds.

**Theorem 6**. *Let* $f_1,\ f_1^*$ *be HFIF for the extended data sets* $P$ *and* $P_{z^*}$, *respectively, and* $\max\{\omega,\ \tilde{\omega}\} < \dfrac{|I_{\min}|}{2|I_{\max}|}$. *Then*

$$\| f_1 - f_1^* \|_\infty \le \frac{1+2\omega-\tilde{\omega}}{1-\omega-\tilde{\omega}} \max_i |z_i - z_i^*|.$$

Let us denote by $f_1^*$ the HVFIF constructed using the construction in section 2 for an extended data set $P^* = \{(x_i^*,\ y_i^*,\ z_i^*) \in \mathrm{R}^3;\ i = 0, 1, \cdots, n\}$. Then from the theorem 4, 5 and 6, we have the following theorem.

**Theorem 7.** *Let* $f_1,\ f_1^*$ *be HFIF for the extended data sets* $P$ *and* $P^*$, *respectively, and* $\max\{\omega,\ \tilde{\omega}\} < \dfrac{|I_{\min}|}{2|I_{\max}|}$. *Then*

$$\| f_1 - f_1^* \|_\infty \le \frac{[(1-\tilde{\omega})(L_1 + L_q) + \omega(L_2 + L_{\tilde{q}})] \max_i |x_i - x_i^*|^\tau + (1+2\omega-\tilde{\omega})(\max_i |y_i - y_i^*| + \max_i |z_i - z_i^*|)}{1-\omega-\tilde{\omega}}.$$

## 4. Box-counting dimension of the HVFIF

In this section, we get a lower and upper bounds for the box-counting dimension of the graph of HVFIF in the case where the data set is $P_0=\{(x_0+\frac{x_n-x_0}{n}i,\, y_i)\in \mathrm{R}^2;\, i=0,1,\cdots,n\}$, the extended data set is $P=\{(x_0+\frac{x_n-x_0}{n}i,\, y_i,\, z_i)\in \mathrm{R}^3;\ \ i=0,1,\cdots,n\}$ and $s_i(x)s_i'(x)\ge 0,\ \tilde{s}_i(x)\tilde{s}_i'(x)\ge 0$, $x\in I,\ i=1,2,\cdots n$ in (2) and (3) in section 2.

Let us denote the graph of the HVFIF $f_1$ by $Gr(f_1)$. As usual, the box-counting dimension of the set *A* by $\dim_B A$ is defined by

$$\dim_B A=\lim_{\delta\to 0}\frac{\log N_\delta(A)}{-\log\delta}\ \text{(if this limit exists)},$$

where $N_\delta(A)$ is any of the following(see [7]):

(i) the smallest number of closed balls of radius $\delta$ that cover the set *A*;
(ii) the smallest number of cubes of side $\delta$ that cover the set *A*;
(iii) the number of $\delta$-mesh cubes that intersect the set *A*;
(iv) the smallest number of sets of diameter at most $\delta$ that cover the set *A*;
(v) the largest number of disjoint balls of radius $\delta$ with centers in the set *A*.

**Lemma 2 (Perron-Frobenius Theorem).** ([16]) *Let $A\ge 0$ be an irreducible square matrix. Then we have the following two statements.*

*(1) The spectral radius $\rho(A)$ of $A$ is an eigenvalue of $A$ and it has strictly positive eigenvector $y$ (i.e., $y_i>0$ for all $i$).*

*(2) $\rho(A)$ increases if any element of $A$ increases.*

Let us denote as follows:

$$\underline{\omega}_k=\min_{x\in I}\{|s_k(x)|,|s_k'(x)|\},\ \ \underline{\tilde{\omega}}_k=\min_{x\in I}\{|\tilde{s}_k(x)|,|\tilde{s}_k'(x)|\},$$

$$\underline{\lambda}=\sum_{i=1}^{n}(\underline{\omega}_i+\underline{\tilde{\omega}}_i),\ \ \overline{\lambda}=\sum_{i=1}^{n}(\overline{\omega}_i+\overline{\tilde{\omega}}_i),$$

where $\overline{\tilde{\omega}}_k$ and $\overline{\omega}_k$ are given in section 3.

For a set $D\subset R$ and a function $f$ defined on $D$, denote as follows:

$$R_f[D]=\sup\{|f(x_2)-f(x_1)|: x_1, x_2\in D\}.$$

**Theorem 8.** *Let $f_1(x)$ be the HVFIF in the Theorem 3. Suppose that there exist three interpolation points $(x_{\alpha_1},y_{\alpha_1})$, $(x_{\alpha_2},y_{\alpha_2})$, $(x_{\alpha_3},y_{\alpha_3})\in P_0$ $(x_{\alpha_1}<x_{\alpha_2}<x_{\alpha_3})$ which are not collinear and tha+t take $z_{\alpha_1}$, $z_{\alpha_2}$ and $z_{\alpha_3}$ such that $(y_{\alpha_i}-y_{\alpha_j})\,(z_{\alpha_i}-z_{\alpha_j})>0$, $i,j=1,2,3, i\ne j$ and three points $(x_{\alpha_1},z_{\alpha_1})$, $(x_{\alpha_2},z_{\alpha_2})$ and $(x_{\alpha_3},z_{\alpha_3})$ are not collinear. Then the box-counting dimension of the graph of $f_1(x)$ is as follows:*

*(a) If $\underline{\lambda}>1$, then $1+\log_n\underline{\lambda}\le\dim_B Gr(f_1)\le 1+\log_n\overline{\lambda}$,*

*(b) If $\overline{\lambda}<1$, then $\dim_B Gr(f_1)=1$.*

**Proof.** Firstly, we prove (*a*). We denote the y-axis vertical distance from the point $(x_{\alpha_2},\ y_{\alpha_2})$ to the line through the points $(x_{\alpha_1},\ y_{\alpha_1})$ and $(x_{\alpha_3},\ y_{\alpha_3})$ by $H$ and the z-axis vertical distance from the point $(x_{\alpha_2},\ z_{\alpha_2})$ to the line through the points $(x_{\alpha_1},\ z_{\alpha_1})$ and $(x_{\alpha_3},\ z_{\alpha_3})$ by $h$. Then we have $H \cdot h > 0$.

We apply $W_i$ to the interpolation points in $I$ one time. Then we have

$$\begin{aligned} &F_i^1(x,y,z) - F_i^1(x,y',z') = \\ &\quad = (s_i(L_i(x))y + s_i'(L_i(x))z + q_i(x)) - (s_i(L_i(x))y' + s_i'(L_i(x))z' + q_i(x)) = \\ &\quad = s_i(L_i(x))(y-y') + s_i'(L_i(x))(z-z'), \\ &F_i^2(x,y,z) - F_i^2(x,y',z') = \tilde{s}_i(L_i(x))(y-y') + \tilde{s}_i'(L_i(x))(z-z'). \end{aligned} \tag{15}$$

Since $Gr(f_1)$ is the graph of a continuous function defined on $I$, the smallest number of $\varepsilon_r$-mesh squares that cover $I_i \times R \cap Gr(f_1)$ is greater than the smallest number of $\varepsilon_r$- mesh squares necessary to cover the vertical line whose length is $\underline{\omega}_i(H+h)$ and less than the smallest number of $\varepsilon_r$-mesh squares necessary to cover the rectangle $I_i \times R_{f_1}[I_i]$.

By (5), (7) and (8) in section 3, we can easily prove that

$$R_{f_1}[I_i] \le \bar{\omega}_i(R_{f_1}[I] + R_{f_2}[I]) + (L_{s_i} \cdot \overline{f_1} + L_{s_i'} \cdot \overline{f_2} + nL_{q_i}) \cdot \frac{1}{n} \ ,$$

$$R_{f_2}[I_i] \le \bar{\tilde{\omega}}_i(R_{f_1}[I] + R_{f_2}[I]) + (L_{\tilde{s}_i} \cdot \overline{f_1} + L_{\tilde{s}_i'} \cdot \overline{f_2} + nL_{\tilde{q}_i}) \cdot \frac{1}{n} \ .$$

Then by (9) and (11) in section 3, we have

$$R_{f_1}[I_i] \le \bar{\omega}_i(R_{f_1}[I] + R_{f_2}[I]) + \frac{1}{n}M,\ \ R_{f_2}[I_i] \le \bar{\tilde{\omega}}_i(R_{f_1}[I] + R_{f_2}[I]) + \frac{1}{n}M \ . \tag{16}$$

Therefore, by (16), we get

$$\sum_{i=1}^{n}\left[\frac{\underline{\omega}_i(H+h)}{\varepsilon_r}\right] \le N(\varepsilon_r) \le \sum_{i=1}^{n}\left(\left[\frac{\bar{\omega}_i(R_{f_1}[I]+R_{f_2}[I]) + M/n}{\varepsilon_r}\right]+1\right)\left(\left[\frac{1}{n\varepsilon_r}\right]+1\right),$$

$$\sum_{i=1}^{n}\left(\frac{\underline{\omega}_i(H+h)}{\varepsilon_r}-1\right) \le N(\varepsilon_r) \le \sum_{i=1}^{n}\left(\frac{\bar{\omega}_i(R_{f_1}[I]+R_{f_2}[I]) + M/n}{\varepsilon_r}+1\right)\left(\frac{1}{n\varepsilon_r}+1\right),$$

$$\frac{\Phi(\mathrm{H}_1)}{\varepsilon_r} - n \le N(\varepsilon_r) \le \left(\frac{\Phi(\mathrm{U}_1)}{\varepsilon_r}+n\right)\left(\frac{1}{n\varepsilon_r}+1\right),$$

where $\Phi(a) = \sum_{i=1}^{n} a_i \ \ (a = (a_1, a_2, \cdots, a_n)^T)$, $\mathrm{H}_1 = (\underline{\omega}_1(H+h), \underline{\omega}_2(H+h), \cdots, \underline{\omega}_n(H+h))^T$,

$$\begin{aligned} \mathrm{U}_1 = \Bigg( &(\bar{\omega}_1 + \bar{\tilde{\omega}}_1)(R_{f_1}[I] + R_{f_2}[I]) + \frac{2M}{n}, (\bar{\omega}_2 + \bar{\tilde{\omega}}_2)(R_{f_1}[I] + R_{f_2}[I]) + \frac{2M}{n}, \cdots, \\ &(\bar{\omega}_n + \bar{\tilde{\omega}}_n)(R_{f_1}[I] + R_{f_2}[I]) + \frac{2M}{n}\Bigg)^T . \end{aligned}$$

Let us apply $W_j$ once agian. Then we have *n* subintervals in every $I_i$. By (15), we get

$$\begin{aligned} F_{ij}^1(x,y,z) - F_{ij}^1(x,y',z') &= [s_j(L_j \circ L_i(x)) f_1^i(x,y,z) + s_j'(L_j \circ L_i(x)) f_2^i(x,y,z) + q_i(L_i(x))] - \\ &\quad - [s_j(L_j \circ L_i(x)) f_1^i(x,y',z') + s_j'(L_j \circ L_i(x)) f_2^i(x,y',z') + q_i(L_i(x))] = \\ &= [s_j(L_j \circ L_i(x)) s_i(L_i(x)) + s_j'(L_j \circ L_i(x)) \tilde{s}_i(L_i(x))](y-y') + \\ &\quad + [s_j(L_j \circ L_i(x)) s_i'(L_i(x)) + s_j'(L_j \circ L_i(x)) \tilde{s}_i'(L_i(x))](z-z'), \end{aligned}$$

$$F_{ij}^2(x,y,z)-F_{ij}^2(x,y',z')=[\tilde{s}_j(L_j\circ L_i(x))s_i(L_i(x))+\tilde{s}'_j(L_j\circ L_i(x))\tilde{s}_i(L_i(x))](y-y')+$$
$$+[s'_j(L_j\circ L_i(x))s'_i(L_i(x))+\tilde{s}'_j(L_j\circ L_i(x))\tilde{s}'_i(L_i(x))](z-z').$$

By (16), we get

$$R_{f_1}[I_{ij}]\le \bar{\omega}_i(R_{f_1}[I_i]+R_{f_2}[I_i])+\frac{1}{n^2}M \quad ,R_{f_2}[I_{ij}]\le \bar{\tilde{\omega}}_i(R_{f_1}[I_i]+R_{f_2}[I_j])+\frac{1}{n^2}M \quad . \tag{17}$$

Since the smallest number of $\varepsilon_r$-mesh squares that cover $I_{ij}\times R\cap Gr(f_1)$ is greater than the smallest number of $\varepsilon_r$- mesh squares necessary to cover the vertical line whose length is $\underline{\omega}_j(\underline{\omega}_i+\underline{\tilde{\omega}}_i)(H+h)$ and less than the smallest number of $\varepsilon_r$-mesh squares necessary to cover the rectangle $I_{ij}\times R_{f_1}[I_{ij}]$, by (17), we have

$$\sum_{i=1}^{n}\sum_{j=1}^{n}\left[\frac{\underline{\omega}_j(\underline{\omega}_i+\underline{\tilde{\omega}}_i)(H+h)}{\varepsilon_r}\right]\le N(\varepsilon_r)\le\sum_{i=1}^{n}\sum_{j=1}^{n}\left(\left[\frac{\bar{\omega}_i(R_{f_1}[I_i]+R_{f_2}[I_i])+M/n^2}{\varepsilon_r}\right]+1\right)\left(\left[\frac{1}{n^2\varepsilon_r}\right]+1\right)$$

$$\sum_{i=1}^{n}\sum_{j=1}^{n}\left(\frac{\underline{\omega}_j(\underline{\omega}_i+\underline{\tilde{\omega}}_i)(H+h)}{\varepsilon_r}-1\right)\le N(\varepsilon_r)\le\sum_{i=1}^{n}\sum_{j=1}^{n}\left(\frac{\bar{\omega}_i(R_{f_1}[I_i]+R_{f_2}[I_i])+M/n^2}{\varepsilon_r}+1\right)\left(\frac{1}{n^2\varepsilon_r}+1\right) \tag{18}$$

Let us denote as follows:

$\underline{S}=\mathrm{diag}(\underline{\omega}_1+\underline{\tilde{\omega}}_1,\ \underline{\omega}_2+\underline{\tilde{\omega}}_2,\ \cdots,\ \underline{\omega}_n+\underline{\tilde{\omega}}_n)$, $\bar{S}=\mathrm{diag}(\bar{\omega}_1+\bar{\tilde{\omega}}_1,\ \bar{\omega}_2+\bar{\tilde{\omega}}_2,\ \cdots,\ \bar{\omega}_n+\bar{\tilde{\omega}}_n)$ and $C$ is $n\times n$ matrix whose entries are all 1.

Then we can rewrite (18) as follows:

$$\frac{\Phi(\mathrm{H}_2)}{\varepsilon_r}-n^2\le N(\varepsilon_r)\le\left(\frac{\Phi(\mathrm{U}_2)}{\varepsilon_r}+n^2\right)\left(\frac{1}{n^2\varepsilon_r}+1\right),$$

where $\mathrm{H}_2=\underline{S}C\mathrm{H}_1$ and $\mathrm{U}_2=\bar{S}C\mathrm{U}_1+2M\ I/n$.

Suppose that

$$\varepsilon_r<\frac{1}{n^k}\le n\varepsilon_r\,. \tag{19}$$

If we apply $W_i$ to $I$ $k$ times, we have

$$\frac{\Phi(\mathrm{H}_k)}{\varepsilon_r}-n^k\le N(\varepsilon_r)\le\left(\frac{\Phi(\mathrm{U}_k)}{\varepsilon_r}+n^k\right)\left(\frac{1}{n^k\varepsilon_r}+1\right), \tag{20}$$

where $\mathrm{H}_k=\underline{S}C\mathrm{H}_{k-1}$ and $\mathrm{U}_k=\bar{S}C\mathrm{U}_{k-1}+2M\ I/n$. Then

$$\mathrm{H}_k=\underline{S}C\mathrm{H}_{k-1}=(\underline{S}C)^2\mathrm{H}_{k-2}=\cdots=(\underline{S}C)^{k-1}\mathrm{H}_1, \tag{21}$$

$$\mathrm{U}_k=\bar{S}C\mathrm{U}_{k-1}+\frac{2M}{n}I=(\bar{S}C)^2\mathrm{U}_{k-2}+\frac{2M}{n}\bar{S}CI+\frac{2M}{n}I=\cdots=$$
$$=(\bar{S}C)^{k-1}\mathrm{U}_1+\frac{2M}{n}(\bar{S}C)^{k-2}I+\frac{2M}{n}(\bar{S}C)^{k-3}I+\cdots+\frac{2M}{n}I\,. \tag{22}$$

Since $\underline{S}C$ and $\bar{S}C$ are non-negative irreducible matrices, by Perron-Frobenius theorem, there exist strictly positive eigenvectors $\underline{\mathrm{e}}$ and $\bar{\mathrm{e}}$ of $\underline{S}C$ and $\bar{S}C$ which correspond to eigenvalues $\underline{\lambda}=\sum_{i=1}^{n}(\underline{\omega}_i+\underline{\tilde{\omega}}_i)$ and $\bar{\lambda}=\sum_{i=1}^{n}(\bar{\omega}_i+\bar{\tilde{\omega}}_i)$, respectively, such that

$$\underline{\mathrm{e}}\le\mathrm{H}_1\,,\quad \bar{\mathrm{e}}\ge\mathrm{U}_1\,,\quad \bar{\mathrm{e}}\ge\frac{2M}{n}I\,.$$

Then by (19), (20) and (21), we get

$$N(\varepsilon_r) \ge \frac{\Phi(\mathrm{H}_k)}{\varepsilon_r} - n^k = \frac{\Phi((\underline{S}C)^{k-1}\mathrm{H}_1)}{\varepsilon_r} - n^k \ge \frac{\Phi((\underline{S}C)^{k-1}\underline{\mathrm{e}})}{\varepsilon_r} - n^k =$$

$$= \frac{\underline{\lambda}^{k-1}\Phi(\underline{\mathrm{e}})}{\varepsilon_r} - n^k \ge \frac{\underline{\lambda}^{k-1}}{\varepsilon_r}\left(\Phi(\underline{\mathrm{e}}) - \frac{1}{\underline{\lambda}^{k-1}}\right)$$

and since $\underline{\lambda} > 1$, we have

$$\frac{\log N(\varepsilon_r)}{-\log \varepsilon_r} \ge 1 - \frac{(k-1)\log \underline{\lambda}}{\log \varepsilon_r} - \frac{\log(\Phi(\underline{\mathrm{e}}) - 1/\underline{\lambda}^{k-1})}{\log \varepsilon_r} \ge 1 + \frac{(k-1)\log \underline{\lambda}}{k \log n} - \frac{\log(\Phi(\underline{\mathrm{e}}) - 1/\underline{\lambda}^{k-1})}{\log \varepsilon_r},$$

$$\lim_{\varepsilon_r \to 0} \frac{\log N(\varepsilon_r)}{-\log \varepsilon_r} \ge \lim_{\varepsilon_r \to 0}\left(1 + \frac{(k-1)\log \underline{\lambda}}{k \log n} - \frac{\log(\Phi(\underline{\mathrm{e}}) - 1/\underline{\lambda}^{k-1})}{\log \varepsilon_r}\right) = 1 + \frac{\log \underline{\lambda}}{\log n} = 1 + \log_n \underline{\lambda}. \quad (23)$$

Since $\underline{\lambda} > 1$,we get $\overline{\lambda} \ge \underline{\lambda} > 1$. Then by (19), (20) and (22), we obtain

$$N(\varepsilon_r) \le \left(\frac{\Phi(\mathrm{U}_k)}{\varepsilon_r} + n^k\right)\left(\frac{1}{n^k \varepsilon_r} + 1\right) =$$

$$= \left[\frac{\Phi\left((\overline{S}C)^{k-1}\mathrm{U}_1 + \frac{2M}{n}(\overline{S}C)^{k-2} I + \frac{2M}{n}(\overline{S}C)^{k-3} I + \cdots + \frac{2M}{n} I\right)}{\varepsilon_r} + n^k\right]\left(\frac{1}{n^k \varepsilon_r} + 1\right)$$

$$\le \left[\frac{\Phi\left((\overline{S}C)^{k-1}\overline{\mathrm{e}} + (\overline{S}C)^{k-2}\overline{\mathrm{e}} + (\overline{S}C)^{k-3}\overline{\mathrm{e}} + \cdots + \overline{\mathrm{e}}\right)}{\varepsilon_r} + n^k\right]\left(\frac{1}{n^k \varepsilon_r} + 1\right) =$$

$$= \left(\frac{(\overline{\lambda}^{k-1} + \overline{\lambda}^{k-2} + \cdots + 1)\Phi(\overline{\mathrm{e}})}{\varepsilon_r} + n^k\right)\left(\frac{1}{n^k \varepsilon_r} + 1\right) = \quad (24)$$

$$= \left(\frac{\Phi(\overline{\mathrm{e}})}{\varepsilon_r}\frac{\overline{\lambda}^{k-1}(1 - 1/\overline{\lambda}^k)}{1 - 1/\overline{\lambda}} + n^k\right)\left(\frac{1}{n^k \varepsilon_r} + 1\right) \ge$$

$$\ge \left(\frac{\Phi(\overline{\mathrm{e}})}{\varepsilon_r}\frac{\overline{\lambda}^{k-1}(1 - 1/\overline{\lambda}^k)}{1 - 1/\overline{\lambda}} + \frac{1}{\varepsilon_r}\right)(n+1) = (n+1)\frac{\overline{\lambda}^{k-1}}{\varepsilon_r}\left(\Phi(\overline{\mathrm{e}})\frac{1 - 1/\overline{\lambda}^k}{1 - 1/\overline{\lambda}} + \frac{1}{\overline{\lambda}^{k-1}}\right) = \frac{\overline{\lambda}^{k-1}}{\varepsilon_r}\gamma,$$

where $\gamma = (n+1)\left(\Phi(\overline{\mathrm{e}})\frac{1 - 1/\overline{\lambda}^k}{1 - 1/\overline{\lambda}} + \frac{1}{\overline{\lambda}^{k-1}}\right) > 0$.

Therefore we have

$$\frac{\log N(\varepsilon_r)}{-\log \varepsilon_r} \le 1 - \frac{(k-1)\log \overline{\lambda}}{\log \varepsilon_r} - \frac{\log \gamma}{\log \varepsilon_r} \le 1 + \frac{(k-1)\log \overline{\lambda}}{k \log n} - \frac{\log \gamma}{\log \varepsilon_r},$$

$$\lim_{\varepsilon_r \to 0} \frac{\log N(\varepsilon_r)}{-\log \varepsilon_r} \le \lim_{\varepsilon_r \to 0}\left(1 + \frac{(k-1)\log \overline{\lambda}}{k \log n} - \frac{\log \gamma}{\log \varepsilon_r}\right) = 1 + \frac{\log \overline{\lambda}}{\log n} = 1 + \log_n \overline{\lambda}. \quad (25)$$

By (23) and (25), we get

$$1 + \log_n \underline{\lambda} \le \dim_B Gr(f_1) \le 1 + \log_n \overline{\lambda}.$$

Secondly, we prove (b). Since $\overline{\lambda} < 1$, by (24), we have

$$N(\varepsilon_r) \le \left(\frac{(\overline{\lambda}^{k-1} + \overline{\lambda}^{k-2} + \cdots + 1)\Phi(\overline{\mathrm{e}})}{\varepsilon_r} + n^k\right)\left(\frac{1}{n^k \varepsilon_r} + 1\right) \le \left(\frac{\Phi(\overline{\mathrm{e}})}{\varepsilon_r}\frac{1 - \overline{\lambda}^k}{1 - \overline{\lambda}} + \frac{1}{\varepsilon_r}\right)(n+1) = \frac{\beta}{\varepsilon_r},$$

where $\beta = \left( \frac{1-\bar{\lambda}^k}{1-\bar{\lambda}} \Phi(\bar{\mathrm{e}}) + 1 \right)(n+1) > 0$.

Then we get

$$\frac{\log N(\varepsilon_r)}{-\log \varepsilon_r} \le 1 - \frac{\log \beta}{\log \varepsilon_r}$$

and

$$\dim_B Gr(f_1) = \lim_{\varepsilon_r \to 0} \frac{\log N(\varepsilon_r)}{-\log \varepsilon_r} \le \lim_{\varepsilon_r \to 0} \left( 1 - \frac{\log \beta}{\log \varepsilon_r} \right) = 1.$$

Since $\dim_B Gr(f_1) \ge 1$, we get $\dim_B Gr(f_1) = 1$. □

## 5. Box-counting dimension of the HVBFIF

As mentioned in section 1, the construction and smoothness of HVFIFs are similar to ones of HVBFIFs in [21, 22]. Therefore, we present only the result for the box-counting dimension of HVBFIFs.

5.1. Construction of HVBFIFs. (see [21])

Let $P_0$ denote the following data set in $\mathrm{R}^3$:

$$P_0 = \{(x_i, y_j, z_{ij}) \in \mathrm{R}^3; i = 0,1,\cdots,n,\ j = 0,1,\cdots,m\}\ (x_0 < x_1 < \cdots < x_n,\ y_0 < y_1 < \cdots < y_m)$$

We extend the data set $P_0$ to a data set $P$ in $\mathrm{R}^4$ as follows:

$$P = \{(x_i, y_j, z_{ij}, t_{ij}) \in \mathrm{R}^4; i = 0,1,\cdots,n,\ j = 0,1,\cdots,m\}\ (x_0 < x_1 < \cdots < x_n,\ y_0 < y_1 < \cdots < y_m)$$

where $t_{ij}$ are parameters. Let us denote $\vec{x}_{ij} = (x_i, y_j)$, $\vec{z}_{ij} = (z_{ij}, t_{ij})$, $N_{nm} = \{1,\cdots,n\} \times \{1,\cdots,m\}$, $I_{x_i} = [x_{i-1}, x_i]$, $I_{y_j} = [y_{j-1}, y_j]$, $I_x = [x_0, x_n]$, $I_y = [y_0, y_m]$, $E = [x_0, x_n] \times [y_0, y_m]$ and $E_{ij} = I_{x_i} \times I_{y_j}$.

We define mappings $L_{x_i} : I_x \to I_{x_i}$, $L_{y_j} : I_y \to I_{y_j}$ for $(i, j) \in N_{nm}$ to be contraction homeomorphisms satisfying the following conditions:

$$L_{x_i}(\{x_0, x_n\}) = \{x_{i-1}, x_i\},\quad L_{y_j}(\{y_0, y_m\}) = \{y_{j-1}, y_j\}.$$

Next, we define transformations $\vec{L}_{ij} : E \to E_{ij}$ by $\vec{L}_{ij}(\vec{x}) = (L_{x_i}(x),\ L_{y_j}(y))$. Then, $\vec{L}_{ij}$ maps the end points of $E$ to the end points of $E_{ij}$, i.e.

$$\vec{L}_{ij}(\vec{x}_{\alpha\beta}) = \vec{x}_{ab},\ a \in \{i-1, i\},\ b \in \{j-1, j\},\ \alpha \in \{0, n\},\ \beta \in \{0, m\}.$$

We define mappings $\vec{F}_{ij} : E \times \mathrm{R}^2 \to \mathrm{R}^2$, $i = 1,\cdots,n$, $j = 1,\cdots,m$ as follows:

$$\vec{F}_{ij}(\vec{x},\ \vec{z}) = \begin{pmatrix} s_{ij}(L_{ij}(\vec{x}))z + s'_{ij}(L_{ij}(\vec{x}))t + q_{ij}(\vec{x}) \\ \tilde{s}_{ij}(L_{ij}(\vec{x}))z + \tilde{s}'_{ij}(L_{ij}(\vec{x}))t + \tilde{q}_{ij}(\vec{x}) \end{pmatrix} = \begin{pmatrix} s_{ij}(L_{ij}(\vec{x})) & s'_{ij}(L_{ij}(\vec{x})) \\ \tilde{s}_{ij}(L_{ij}(\vec{x})) & \tilde{s}'_{ij}(L_{ij}(\vec{x})) \end{pmatrix} \begin{pmatrix} z \\ t \end{pmatrix} + \begin{pmatrix} q_{ij}(\vec{x}) \\ \tilde{q}_{ij}(\vec{x}) \end{pmatrix} \tag{26}$$

where $\vec{x} = (x, y)$, $\vec{z} = (z, t)$ and $s_{ij}, s'_{ij}, \tilde{s}_{ij}, \tilde{s}'_{ij} : E_{ij} \to R$ are arbitrary Lipschitz functions whose absolute values are less than 1 and $q_{ij}$, $\tilde{q}_{ij} : E \to R$ are Lipschitz functions satisfying the following conditions: If $\alpha \in \{0, n\}$, $\beta \in \{0, m\}$, $a \in \{i-1, i\}$, $b \in \{j-1, j\}$, $L_{x_i}(x_\alpha) = x_a$, $L_{y_j}(y_\beta) = y_b$, then $\vec{F}_{ij}(\vec{x}_{\alpha\beta}, \ \vec{z}_{\alpha\beta}) = \vec{z}_{ab}$.

We assume that $D \subset R^2$ is an enough large bounded domain containing $\vec{z}_{ij}$, $i = 1, \cdots, n$, $j = 1, \cdots, m$. We define transformations $\vec{W}_{ij} : E \times D \to E_{ij} \times \mathrm{R}^2$ as follows:

$$\vec{W}_{ij}(\vec{x}, \vec{z}) = (\vec{L}_{ij}(\vec{x}), \ \vec{F}_{ij}(\vec{x}, \vec{z})), \ \ i = 1, \cdots, n, \ \ j = 1, \cdots, m.$$

For a function $f$, let us denote $\overline{f} = \sup_x |f(x)|$ and $\overline{S} = \max\{\bar{s}_{ij} + \bar{\tilde{s}}_{ij}, \bar{s}'_{ij} + \bar{\tilde{s}}'_{ij}; i = 1, \cdots, n$, $j = 1, \cdots, m\}$. If $\overline{S} < 1$, then there exists a distance $\rho_\theta$ equivalent to the Euclidean metric on $\mathrm{R}^2$ such that $\vec{W}_{ij}$, $i = 1, \cdots, n$, $j = 1, \cdots, m$ are contraction transformations with respect to the distance $\rho_\theta$. (See Theorem 2.1 in [21].)

Therefore, we have an hyperbolic iterated function system $\{R^4; \vec{W}_{ij}, \ i = 1, \cdots, n, \ \ j = 1, \cdots, m\}$ corresponding to the extended data set $P$. Then there exists a continuous interpolation function $\vec{f}$ of the extended data set $P$ such that the graph of $\vec{f}$ is the attractor of the IFS. (See Theorem 3.1 in [21].)

We denote the vector valued function $\vec{f} : E \to R^2$ by $\vec{f} = (f_1(x, y), \ f_2(x, y))$, where $f_1 : E \to R$ interpolates the given data set $P_0$ and the function $f_2(x, y)$ interpolates the set $\{(x_i, y_j, t_{ij}) = (\vec{x}_{ij}, t_{ij}) \in \mathrm{R}^3; \quad i = 0, 1, \cdots, n, \quad j = 0, 1, \cdots, m\}$. Then, we have

$$\vec{f}(x, \ y) = \vec{F}_{ij}(\vec{L}_{ij}^{-1}(x, y), \ \vec{f}(\vec{L}_{ij}^{-1}(x, y))) = \vec{F}_{ij}(\vec{L}_{ij}^{-1}(x, y), \ f_1(\vec{L}_{ij}^{-1}(x, y)), \ f_2(\vec{L}_{ij}^{-1}(x, y))), (x, y) \in E_{ij},$$

$$f_1(x, y) = s_{ij}(x, y) f_1(\vec{L}_{ij}^{-1}(x, y)) + s'_{ij}(x, y) f_2(\vec{L}_{ij}^{-1}(x, y)) + q_{ij}(\vec{L}_{ij}^{-1}(x, y)),$$

$$f_2(x, y) = \tilde{s}_{ij}(x, y) f_1(\vec{L}_{ij}^{-1}(x, y)) + \tilde{s}'_{ij}(x, y) f_2(\vec{L}_{ij}^{-1}(x, y)) + \tilde{q}_{ij}(\vec{L}_{ij}^{-1}(x, y)).$$

5.2. Box-counting dimension of HVBFIFs

Let the data set be $P_0=\left\{\left(x_0+\frac{x_n-x_0}{n}i,\, y_0+\frac{y_n-y_0}{n}j,\, z_{ij}\right)\in \mathrm{R}^3; i,\, j=0,1,\cdots,n\right\}$ and the extended data set $P=\left\{\left(x_0+\frac{x_n-x_0}{n}i,\, y_0+\frac{y_n-y_0}{n}j,\, z_{ij},t_{ij}\right)\in \mathrm{R}^4; i,\, j=0,1,\cdots,n\right\}$ and $s_{ij}(x,y)s'_{ij}(x,y)\ge 0$, $\tilde{s}_{ij}(x,y)\tilde{s}'_{ij}(x,y)\ge 0$, $(x,y)\in E,\ i,\, j=1,2,\cdots,n$ in (26).

Let us denote as follows:

$$P_{0x_\alpha}=\left\{\left(x_0+\frac{x_n-x_0}{n}\alpha,\, y_0+\frac{y_n-y_0}{n}j,\, z_{\alpha j}\right)\in \mathrm{R}^3;\ j=0,1,\cdots,n\right\},\ \alpha=0,1,...,n,$$

$$P_{0y_\beta}=\left\{\left(x_0+\frac{x_n-x_0}{n}i,\, y_0+\frac{y_n-y_0}{n}\beta,\, z_{i\beta}\right)\in \mathrm{R}^3; i=0,1,\cdots,n\right\},\ \beta=0,1,...,n,$$

$$\underline{\omega}_{ij}=\min_{(x,y)\in E_{ij}}\{|s_{ij}(x,y)|,|s'_{ij}(x,y)|\},\quad \underline{\tilde{\omega}}_{ij}=\min_{(x,y)\in E_{ij}}\{|\tilde{s}_{ij}(x,y)|,|\tilde{s}'_{ij}(x,y)|\},$$

$$\overline{\omega}_{ij}=\max_{(x,y)\in E_{ij}}\{|s_{ij}(x,y)|,|s'_{ij}(x,y)|\},\quad \overline{\tilde{\omega}}_{ij}=\max_{(x,y)\in E_{ij}}\{|\tilde{s}_{ij}(x,y)|,|\tilde{s}'_{ij}(x,y)|\},$$

$$\underline{\lambda}=\sum_{i,j=1}^{n}(\underline{\omega}_{ij}+\underline{\tilde{\omega}}_{ij}),\quad \overline{\lambda}=\sum_{i,j=1}^{n}(\overline{\omega}_{ij}+\overline{\tilde{\omega}}_{ij}).$$

**Theorem 9.** Let $f_1(x,y)$ be the HVBFIFs constructed above. Suppose that there exist three interpolation points $(x_\alpha,y_{j_1},z_{\alpha j_1})$, $(x_\alpha,y_{j_2},z_{\alpha j_2})$, $(x_\alpha,y_{j_3},z_{\alpha j_3})\in P_{0x_\alpha}$ $(y_{j_1}<y_{j_2}<y_{j_3})$ (or $(x_{i_1},y_\beta,z_{i_1\beta})$, $(x_{i_2},y_\beta,z_{i_2\beta}),(x_{i_3},y_\beta,z_{i_3\beta})\in P_{0y_\beta}$ $(x_{i_1}<x_{i_2}<x_{i_3})$) which are not collinear and that take $t_{\alpha j_1}$, $t_{\alpha j_2}$ and $t_{\alpha j_3}$ (or $t_{i_1\beta}$, $t_{i_2\beta}$ and $t_{i_3\beta}$) such that $(z_{\alpha j_k}-z_{\alpha j_l})\ (t_{\alpha j_k}-t_{\alpha j_l})>0$ (or $(z_{i_k\beta}-z_{i_l\beta})\ (t_{i_k\beta}-t_{i_l\beta})>0$), $k,l=1,2,3,k\ne l$ and three points $(x_\alpha,y_{j_1},t_{\alpha j_1})$, $(x_\alpha,y_{j_2},t_{\alpha j_2})$, $(x_\alpha,y_{j_3},t_{\alpha j_3})$ (or $(x_{i_1},y_\beta,t_{i_1\beta})$, $(x_{i_2},y_\beta,t_{i_2\beta})$, $(x_{i_3},y_\beta,t_{i_3\beta})$) are not collinear. Then the box-counting dimension of the graph of $f_1(x,y)$ is as follows:

*(a) If* $\underline{\lambda}>n$, *then* $1+\log_n\underline{\lambda}\le \dim_B Gr(f_1)\le 1+\log_n\overline{\lambda}$,

*(b) If* $\overline{\lambda}\le n$, *then* $\dim_B Gr(f_1)=2$.

## 6. Conclusion

In the present paper, we introduce a construction of hidden variable fractal interpolation functions (HVFIFs) in $R^2$ by IFSs with four function contractivity factors and analyze their smoothness and stability, which are similar to ones in [21, 22]. We obtain the lower and upper bounds of box-counting dimension of the constructed HVFIFs and the HVBFIFs introduced in [21].

The results of the HVFIFs and HVBFIFs become theoretical foundation for their practical applications such as approximation of functions and curves, interpolation, image processing, data fitting and computer graphics, etc.